\documentclass[12pt]{amsart}
\usepackage{amsmath, amscd}
\usepackage{amsfonts}\usepackage{ulem}
\usepackage{amssymb}
\usepackage{xcolor}
\usepackage{marginnote}

\newcommand{\de}{\partial}
\newcommand{\R}{\mathbb R}

\newcommand{\al}{\alpha}

\newcommand{\C}{\mathbb C}

\newcommand{\B}{\mathbb B}

\newcommand{\D}{\mathbb D}

\newcommand{\N}{\mathbb N}

\def\vp{\varphi}

\def\Re{{\sf Re}\,}

\newtheorem{theorem}{Theorem}[section]

\newtheorem{proposition}[theorem]{Proposition}

\newtheorem{question}[theorem]{Question}

\theoremstyle{definition}
\newtheorem{definition}[theorem]{Definition}
\newtheorem{example}[theorem]{Example}

\theoremstyle{remark}
\newtheorem{remark}[theorem]{Remark}
\numberwithin{equation}{section}

\numberwithin{equation}{section}

\begin{document}
\title[Support points and the Bieberbach conjecture]{Support points and the Bieberbach conjecture in higher dimension}
\author[F. Bracci]{Filippo Bracci$^\dag$} \address{F. Bracci: Dipartimento Di Matematica\\
Universit\`{a} di Roma \textquotedblleft Tor Vergata\textquotedblright\ \\
Via Della Ricerca Scientifica 1, 00133 \\
Roma, Italy} \email{fbracci@mat.uniroma2.it}
\author[O. Roth]{Oliver Roth} \address{O. Roth: Department of Mathematics,
  University of W\"urzburg, Emil Fischer Stra{\ss}e 40, 97074 W\"urzburg,
  Germany} \email{roth@mathematik.uni-wuerzburg.de}

\thanks{$^\dag$Supported by the ERC grant ``HEVO - Holomorphic Evolution Equations'' n. 277691}
\begin{abstract}
We describe some open questions related to support points in the class $S^0$ and introduce some useful techniques  toward a higher dimensional Bieberbach conjecture. 
\end{abstract}
\maketitle
\tableofcontents

\section{Introduction}

Let  $\B^n\subset \C^n$ denote the unit ball for the standard Hermitian product in $\C^n$, $n\geq 1$. For the sake of simplicity, we consider only the case $n=1$ and $n=2$, but, in fact, all results which we are going to discuss for $\B^2$ hold  in any dimension. We  denote $\D:=\B^1$.

Let $S(\B^n)$ denote the class of univalent maps $f:\B^n \to \C^n$ normalized so that $f(0)=0, df_0={\sf id}$. We consider $S(\B^n)$ as a subspace of the Frech\'et space of holomorphic maps from $\B^n$ to $\C^n$ with the topology of uniform convergence on compacta. 

For $n=1$, the set $S(\D)$ is compact (see,  e.g., \cite{P}). For $n>1$, the set $S(\D)$ is not compact: a simple example is given by considering the sequence 
\[
\{(z_1,z_2)\mapsto (z_1+m z_2^2, z_2)\}_{m\in \N},
\]
 which belongs to $S(\B^2)$ but for $m\to \infty$ does not converge. 

Since the set $S(\D)$ is compact, for every continuous linear operator $L: {\sf Hol}(\D, \C)\to \C$ there exists $f\in S(\D)$ such that $\Re L(f)\geq \Re L(g)$ for all $g\in S(\D)$.  If $L$ is not constant, such an $f$ is called a {\sl support point} for $L$.

It is known that every support point in $S(\D)$ is unbounded and is, in fact,  a slit map  (see \cite{Sp}). The most interesting linear functionals to be considered are perhaps those defined as $L_m(f)=b_m$, where $f(z)=\sum_{j\geq 0} b_j z^j \in S(\D)$, $m\in \N$. The Bieberbach conjecture, proved in the 80's by L. de Branges, states that $|b_m|\leq m$ for $f\in S(\D)$ and that equality is reached only by rotations of the Koebe function. 

In higher dimensions, since the class $S(\B^2)$ is not compact, one is forced to consider suitable compact subclasses. Convex maps and starlike maps form compact subclasses, but, for many purposes, these classes are too small. In \cite{GHK} (see also \cite{GK}), it was introduced a compact subclass, denoted by $S^0(\B^2)$ (or simply $S^0$), for which the membership  depends  on the existence of a parametric representation, a condition that is always satisfied in dimension one thanks to the classical Loewner theory (see Section \ref{defS0}). 

 The class $S^0$ is strictly
 contained in $S:=S(\B^2)$, but evidences are that every map
 in  $S$ might be factorized as the composition of an element of $S^0$ and a normalized entire univalent map of $\C^2$ (this is known to be true for univalent maps $f$ on $\B^2$ which extend $C^\infty$ up to the boundary, $f(\B^2)$ is strongly pseudoconvex and $\overline{f(\B^2)}$ is polynomially convex\footnote{in \cite{I}, the result was extended to univalent maps on $\B^2$ which extend $C^1$ up to $\partial \B^2$ and whose image is Runge in $\C^2$, but, unfortunately, there is a gap in the proof};  see \cite{ABW2}). Were this the case, one could somehow split the difficulties in understanding univalent maps on $\B^2$ into two pieces: understanding the compact class $S^0$ and automorphisms/Fatou-Bieberbach maps in $\C^2$. In  light of this, the class $S^0$ seems to be a natural candidate to study in higher dimensions. 

The present note focuses on support points on $S^0$. Our aim is, on the one hand, to state some natural open questions originating in the recent works \cite{BGHK, B, I, R}, and, on the other hand, to develop some new techniques to handle such problems (in particular, slice reduction and decoupling harmonic terms tricks). The paper \cite{GHKK1} contains other open questions in this direction and an extensive bibliography on the subject, to which we refer the reader. Here we mainly focus our attention on those questions which relate the class $S^0$ to the (huge) group of automorphisms of $\C^2$, highlighting the deep differences between dimension $1$ and dimension $2$ (see Section \ref{auto}).

 We also develop the ideas in \cite{B} (see Section \ref{operaz}), which allowed to construct an example of a bounded support point in $S^0$. With these tools in hand, we state a Bieberbach-type conjecture in $S^0$ for coefficients of pure terms in $z_1$ and $z_2$ in the expansion at the origin. 
 
\medskip

We thank the referee for the comments which improved the original manuscript.

\section{The class $S^0$}\label{defS0}
In what follows we denote by $\R^+$ the semigroup of nonnegative real numbers and by $\N$ the semigroup of nonnegative integers.

Let
\[\mathcal M:=\{h\in {\sf Hol}(\B^2, \C^2): h(0)=0, dh_0={\sf id},
\Re \langle h(z), z\rangle >0,  \forall z\in \B^2\setminus\{0\}\},\]
where $\langle \cdot, \cdot \rangle$ denotes the Euclidean inner product in $\C^2$.

The set $\mathcal M$ is compact in ${\sf Hol}(\B^2, \C^2)$ endowed with the topology of uniform convergence on compacta (see \cite{GK}).

\begin{definition}\label{Herglotz}
 A \textit{Herglotz vector field associated with the class $\mathcal M$} on $\B^2$ is a mapping $G:\B^2\times \R^+\to
\C^2$ with the following properties:
\begin{itemize}
\item[(i)] the mapping $G(z,\cdot)$ is measurable on $\R^+$ for all
$z\in \B^2$.
\item[(ii)] $-G(\cdot, t)\in \mathcal M$  for a.e. $t\in [0,+\infty)$.
\end{itemize}
\end{definition}

\begin{remark}
Due to the estimates for the class $\mathcal M$ (see \cite{GK}),  a Herglotz vector field associated with the class $\mathcal M$ on $\B^2$ is an $L^\infty$-Herglotz vector field on $\B^2$ in the sense of \cite{BCD}.
\end{remark}

\begin{definition}\label{regular}
A family $(f_t)_{t\geq 0}$ of  holomorphic mappings from $\B^2$ to $\C^2$ such that $f_t(0)=0$ and $d(f_t)_0=e^t{\sf id}$ for all $t\geq 0$, is called a {\sl normalized regular family} if
\begin{itemize}
\item[(i)]   the mapping  $t\mapsto f_t$ is continuous with respect to the topology  in ${\sf Hol}(\B^2,\C^2)$ induced by the uniform convergence on compacta in $\B^2$,
\item[(ii)] there exists a set of zero measure $N\subset [0,+\infty)$
such that for all $t\in [0,+\infty)\setminus N$ and all $z\in \B^2$
the partial derivative $\frac{\de f_t}{\de t}(z)$ exists and  is holomorphic.
\end{itemize}
\end{definition}

For a given Herglotz vector field $G(z,t)$ associated with the class $\mathcal M$ on $\B^2$, a {\sl normalized solution} to the  {\sl Loewner-Kufarev PDE} associated to $G(z,t)$ consists of a normalized regular family $(f_t)_{t\geq 0}$ such that  the following equation is satisfied for a.e. $t\geq 0$ and for all  $z\in \B^2$
\begin{equation}\label{LoewnerPDE}
\frac{\de f_t}{\de t}(z)=-d(f_t)_z \cdot G(z,t).
\end{equation}

\begin{definition}
A {\sl normalized subordination chain} $(f_t)_{t\geq 0}$ is a family of holomorphic mappings
 $f_t:{\mathbb B^2}\rightarrow {\mathbb C^2}$, such that $f_t(0)=0$, $d(f_t)_0=e^t {\sf id}$ for all $t\geq 0$,
and for every $0\leq s\leq t$ there exists $\vp_{s,t}:\B^2 \to \B^2$ holomorphic such that $f_s=f_t\circ \vp_{s,t}$.
A normalized subordination chain $(f_t)_{t\geq 0}$ is called a {\sl normalized Loewner chain} if for all $t\geq 0$ the mapping $f_t$ is univalent.
\end{definition}

\begin{definition}
A normalized Loewner chain $(f_t)_{t\geq 0}$ on $\B^2$ is called a
{\sl normal Loewner chain} if the family $\{e^{-t}f_t(\cdot)\}_{t\geq 0}$ is normal.
\end{definition}

From \cite[Chapter 8]{GK},  \cite[Prop. 2.6]{ABW} and \cite{GKP}, we have the following:

\begin{theorem}\label{all}
\begin{enumerate}
  \item If $(f_t)_{t\geq 0}$ is a normalized Loewner chain on $\B^2$,
  then it is a normalized solution to a Loewner-Kufarev PDE \eqref{LoewnerPDE}
  for some Herglotz vector field $G(z,t)$ associated with the class $\mathcal M$ in $\B^2$.
  \item Let $G(z,t)$ be a Herglotz vector field associated with
  the class $\mathcal M$ on $\B^2$. Then there exists a unique
  normal Loewner chain $(g_t)_{t\geq 0}$---called the {\sl canonical solution}---which is a
  normalized solution to \eqref{LoewnerPDE}. Moreover, $\bigcup_{t\geq 0}g_t(\B^2)=\C^2$.
  \item If $(f_t)_{t\geq 0}$ is a normalized solution to \eqref{LoewnerPDE},
  then $(f_t)_{t\geq 0}$ is a normalized subordination chain on $\B^2$.
  Moreover, there exists a holomorphic mapping $\Phi: \C^2 \to \bigcup_{t\geq 0} f_t(\B^2)$,
  with $\Phi(0)=0$ and $d\Phi_0={\sf id}$ such that $f_t=\Phi \circ g_t$, where $(g_t)_{t\geq 0}$
  is the canonical solution to \eqref{LoewnerPDE}.
  In particular, $(f_t)_{t\geq 0}$ is a normalized Loewner chain if and only if $\Phi$ is univalent.
\end{enumerate}
\end{theorem}

\begin{remark}
Let $(f_t)_{t\geq 0}$  be a family of holomorphic mappings such that
$f_t(0)=0, d(f_t)_0=e^t{\sf id}$ for all $t\geq 0$. Assume that $(f_t)_{t\geq 0}$
satisfies (i) of Definition \ref{regular} and for all fixed $z\in \B^2$
the mapping $t\mapsto f_t(z)$ is absolutely continuous.
If for all fixed $z\in \B^2$ the family $(f_t)_{t\geq 0}$
satisfies \eqref{LoewnerPDE} for a.e. $t\geq 0$, then it a normalized solution to the
Loewner-Kufarev PDE and thus  is a regular family.
\end{remark}

\begin{definition}
Let $f\in S$. We say that $f$ admits {\sl parametric representation} if
$$f(z)=\lim_{t\to\infty}e^t\varphi(z,t)$$
locally uniformly on $\B^2$, where
$\varphi(z,0)=z$ and
\begin{equation}\label{LLODE}
\frac{\partial \varphi}{\partial t}(z,t)=G(\varphi(z,t),t),\quad \mbox{ a.e. }\, t\geq 0,\quad
\forall\, z\in\B^2,
\end{equation}
for some Herglotz vector field $G$ associated with the class ${\mathcal M}$ on $\B^2$.
\end{definition}

We denote by $S^0$ the set consisting of univalent mappings which admit parametric representation.

The following result is in \cite[Chapter 8]{GK} (see also \cite{GHK}):

\begin{theorem}
\begin{enumerate}
\item A normalized univalent map $f: \B^2\to \C^2$ has parametric representation if and only if
there exists a normal Loewner chain $(f_t)_{t\geq 0}$ on $\B^2$
such that $f_0=f$.
\item The class $S^0$ is compact in the topology of uniform
convergence on compacta.
\end{enumerate}
\end{theorem}

\section{Support points and extreme points}

\begin{definition}
(i)
Let $K$ be a compact subset of ${\sf Hol}(\B^2, \C^2)$ endowed with the
topology of uniform convergence on compacta.
  A mapping $f\in K$ is called a {\sl support point} if there exists a continuous linear operator $L:{\sf Hol}(\B^2, \C^2)\to \C$  not constant on $K$ such that $\max_{g\in K} \Re L(g)=\Re L(f)$. We denote by ${\sf Supp}(K)$ the set of support points of $K$.

(ii)
A mapping $f\in K$ is called an {\sl extreme point}
if $f=tg+(1-t)h$, where $t\in (0,1)$, $g,h\in K$, implies
$f=g=h$.
We denote by ${\sf Ex}(K)$ the set of extreme points of $K$.
\end{definition}

Note that the notion of extreme points is not related to topology, but only on the geometry of the set. If $K$ is a compact subset of ${\sf Hol}(\B^2, \C^2)$ and $a\in K$ is a support point which maximizes the continuous linear operator $L$, then $\mathcal L:=\{b\in {\sf Hol}(\B^2, \C^2): \Re L(b)=\Re L(a)\}$ is a real hyperplane and $\mathcal L\cap K$ contains extreme points. Therefore, for any continuous linear operator $L$ which is not constant on $K$ there exists a point $a\in K$ which is  both a support point (for $L$) and an extreme point for $K$.

In dimension one it is known that all support points for $S^0$ are slit mappings (see \cite{Sp}). In higher dimension, the situation is considerably more complicated.

\begin{proposition}\cite{R}
Let $f\in S^0$ be a support point. Let $G(z,t)$ be a Herglotz vector field
associated with the class $\mathcal M$ which generates a normal Loewner chain $(f_t)$ such that $f_0=f$. Then $G(z,t)$ is a support point of $-\mathcal M$ for a.e. $t\geq 0$.
\end{proposition}

\begin{question} 
Let $f\in S^0$ be a support point. 
\begin{enumerate}
\item Does there exist only one normal Loewner chain $(f_t)$ such that
$f_0=f$\,?
\item Does there exist a
  Herglotz vector field associated with the class $\mathcal M$  which  generates a normal Loewner chain $(f_t)$ with
$f_0=f$ such that $t
  \mapsto G(\cdot,t)$ is continuous and $G(\cdot,t) \in {\sl supp} (-\mathcal{M})$
for all $t \ge 0$\,?
\end{enumerate}
\end{question}

\begin{question}
Let $G(z,t)$ be a Herglotz vector field associated with the class $\mathcal M$ which generates a normal Loewner chain $(f_t)$. 
\begin{enumerate}
  \item If $f_0$ is extreme in $S^0$, is it true that $G(z,t)$ is extreme in $-\mathcal M$ for a.e. $t\geq 0$?
\item If $f_0$ is extreme in $S^0$, is  $G(z,t)$ 
  uniquely determined\,?
\end{enumerate}
\end{question}

\begin{proposition}\cite{Schl}
Let $(f_t)$ be a normal Loewner chain. Then for all $t\geq 0$, $e^{-t}f_t\in S^0$. Moreover, if $f_0$ is a support/extreme point for $S^0$, so is $e^{-t}f_t$ for all $t\geq 0$.
\end{proposition}

\begin{definition}
Let $(f_t)_{t\geq 0}$ be a normalized Loewner chain in $\B^2$ and $G(z,t)$ be the associated Herglotz vector field.
We say that $(f_t)_{t\geq 0}$ is  {\sl exponentially squeezing  in $[T_1,T_2)$,
for $0\leq T_1< T_2\leq +\infty$ (with squeezing ratio  $a\in (0,1))$} 
if for a.e. $t\in [T_1,T_2)$ and for all $z\in \B^2\setminus\{0\}$,
\begin{equation}\label{squeezing}
      \Re \left\langle G(z,t), \frac{z}{\|z\|^2}\right\rangle \leq -a.
\end{equation}
\end{definition}

In \cite{BGHK} it is proved that \eqref{squeezing} is equivalent to: for all $T_1\leq s<t<T_2$, 
\begin{equation}\label{equiv}
\|f_t^{-1}(f_s(z))\|\leq e^{a(s-t)}\|z\|,  \quad \hbox{for all $z\in \B^2$}.
\end{equation}

Hence, if $(f_t)$ is exponentially squeezing in $[T_1,T_2)$, then $f_t$ is bounded for all $t\in [0,T_2)$ and $\overline{f_s(\B^2)}\subset f_t(\B^2)$ for all $T_1\leq s< t<T_2$.

Using the results of \cite{BGHK}, or of \cite{R} for support points, one can prove

\begin{proposition}
Let $(f_t)$ be a normal Loewner chain which is exponentially squeezing in  $[T_1,T_2)$. Then $f_0\not\in {\sf Supp}(S^0)\cup {\sf Ex}(S^0)$.
\end{proposition}

\begin{example}
Let $f\in S^0$. Let $(f_t)$ be one parametric representation of $f$. Let $r\in (0,1)$. Consider $f_{r,t}(z):=r^{-1}f_t(rz)$. Then $(f_{r,t})$ is an exponentially squeezing normal Loewner chain and in particular, $f_r\in S^0\setminus ({\sf Supp}(S^0)\cup {\sf Ex}(S^0))$.
\end{example}

\begin{question}
Let $f\in S^0\setminus ({\sf Supp}(S^0)\cup {\sf Ex}(S^0))$
be a bounded function.  Is it true that $f$ can be embedded into an exponentially squeezing Loewner chain?
\end{question}

\begin{question}
Let $G(z,t)$ be a Herglotz vector field associated with the class $\mathcal M$ which generates a normal Loewner chain $(f_t)$. Assume that
\[
 \limsup_{z\to A}\Re \left\langle G(z,t), \frac{z}{\|z\|^2}\right\rangle \leq -a,
\]
for some $A\subset \de\B^2$ and a.e. $t\geq 0$.
Is it true that $f_0\not\in {\sf Supp}(S^0)\cup {\sf Ex}(S^0)$?
\end{question}

\section{Automorphisms of $\C^2$ and support points}\label{auto}

Let  
\[
{\sf Aut}_0(\C^2):=\{f\in {\sf Aut}(\C^2): f(0)=0, df_0={\sf id}\}.
\]
Given $f\in {\sf Aut}_0(\C^2)$,   for every $r>0$, the map $f^r: \B^2 \to \C^2$ defined by
\[
f^r(z)=\frac{1}{r}f(rz),
\]
is normalized and univalent. For $r<<1$, the image $f^r(\B^2)$ is convex, hence $f^r\in S^0$. For $f\in {\sf Aut}_0(\C^2)$, let
\[
r(f):=\sup\{t>0: f^t\in S^0\}.
\]
Since $S^0$ is compact and $\{f^t\}_{t>0}$ is not normal except for $f={\sf id}$, it follows that for $f\in {\sf Aut}_0(\C^2)\setminus\{{\sf id}\}$
\[
0<r(f)<+\infty, \quad f^{r(f)}\in S^0.
\]

\begin{question}\label{questionauto}
Let $f\in {\sf Aut}_0(\C^2)\setminus\{{\sf id}\}$. Is it true that $f^{r(f)}\in {\sf Supp}(S^0)$?
\end{question}

The previous question has a positive solution in the case $f(z_1,z_2)=(z_1+az_2^2, z_2)$, as we discuss later (or see \cite{B}).

Let
\[
\mathcal A:=\{f \in S^0: \hbox{there exists } \Psi \in {\sf Aut}(\C^2) : \Psi|_{\B^2}=f\}
\] 
Note that in dimension one the  analogue of $\mathcal A$
  contains only the identity mapping.
In higher dimension we have

\begin{theorem}\cite{I} 
$\overline{\mathcal A}=S^0$.
\end{theorem}

Take $f\in S^0$, and expand $f$ as
\[
f(z_1,z_2)=(z_1+\sum_{\al\in \N^2, |\al|\geq 2} b_\al^1 z^\al, z_2+\sum_{\al\in \N^2, |\al|\geq 2} b_\al^2 z^\al).
\]
By a result of F. Forstneri\v{c} \cite{F}, for any $M\in \N$ there exists $g\in {\sf Aut}_0(\C^2)$ such that $f-g=O(\|z\|^{M+1})$ (that is $f, g$ have the same jets up to order $M$). However, such a $g$ does not belong to $S^0$ in general. 
\begin{question}
For which $\al\in \N^2$ is it true that  for any $f\in S^0$ there exists $g\in \mathcal A$ having the same coefficients $b_\al^1$ as $f$?
\end{question}

\section{Coefficient bounds in $\B^2$}

We  use the following notation:  $f\in S^0$, 
\[
f(z_1,z_2)=(z_1+\sum_{\al\in \N^2, |\al|\geq 2} b_\al^1 z^\al, z_2+\sum_{\al\in \N^2, |\al|\geq 2} b_\al^2 z^\al).
\]
If $(f_t)$ is a normal Loewner chain, we denote by $b^j_\al(t)$ the corresponding coefficients of $f_t$.

For $G(z,t)$ a Herglotz vector field associated with the class $\mathcal M$
\[
G(z,t)=(-z_1+\sum_{\al\in \N^2,  |\al|\geq 2} q^1_\al(t) z^\al, -z_2+\sum_{\al\in \N^2,  |\al|\geq 2} q^2_\al(t) z^\al).
\]
For an evolution family $\vp_{s,t}:=f_t^{-1}\circ f_s$,
\[
\vp_{s,t}(z_1,z_2)=(e^{s-t}z_1+\sum_{\al\in \N^2,  |\al|\geq 2} a^1_\al(s,t) z^\al, e^{s-t}z_2+\sum_{\al\in \N^2,  |\al|\geq 2} a^2_\al(s,t) z^\al).
\]
For $f\in S^0$, $\al\in \N^2, |\al|\geq 2$, $j=1,2$, let 
\[
L^j_\al(f):=b^j_\al.
\]
The problem of finding the maximal possible sharp bound for coefficients of mappings in the class $S^0$, consists in fact in finding the support points in $S^0$ for the linear functionals $L^j_\al$.

If $(f_1(z_1,z_2), f_2(z_1,z_2))\in S^0$, then $(f_2(z_2,z_1), f_1(z_2,z_1))\in S^0$. Therefore, it is enough to solve the problem for $L^1_\al$. More generally, if $U$ is a $2\times 2$ unitary matrix, given a normal Loewner chain $(f_t)$, $(U^\ast f_t (Uz))$ is again a normal Loewner chain. This enables us to assume that a given coefficient $b^1_\al>0$, so that, in fact, $\max_{g\in S^0}\Re L^1_\al(g)=\max_{g\in S^0}|b_\al^1(g)|$.

Let $(f_t)$ be a normal Loewner chain, $G(z,t)$ the associated Herglotz vector field and $(\vp_{s,t})$ the associated evolution equation. Expanding the Loewner ODE one gets
\begin{equation}\label{ODE}
\frac{\de a^1_\al(s,t)}{\de t}=-a^1_\al(s,t)+q_\al^1(t)e^{|\al|(s-t)}+R_\al,
\end{equation}
where $R_\al$ is the coefficient of $z^\al$ in the expansion of
\[
\sum_{2\leq |\gamma|\leq |\al|-1}q^1_\gamma(t)(e^{s-t}z_1+\sum_{ 2\leq |\beta|\leq |\al|-1} a^1_\beta(s,t) z^\beta)^{\gamma_1}(e^{s-t}z_1+\sum_{ 2\leq |\beta|\leq |\al|-1} a^2_\beta(s,t) z^\beta)^{\gamma_2}.
\]
Since $f_0=\lim_{t\to\infty} e^{t}\vp_{0,t}$ uniformly on compacta, we have
\begin{equation}\label{b-a}
b^1_\al=\lim_{t\to \infty} e^t a^1_\al(0,t).
\end{equation}
Therefore, in order to get a sharp bound on the coefficients, one should try first to reduce the problem (if possible) to a simple problem involving the least possible number of coefficients of $G$, then find a sharp bound for such coefficients and solve the associated ODE. 

Below, we describe  some methods which can be used  to simplify the problem and then turn to some applications.

\section{Operations in the class $\mathcal M$}\label{operaz}
\subsection{Decoupling harmonic terms} Let $G(z)$ be an autonomous Herglotz vector field associated with the class $\mathcal M$. Then 
\[
\Re \langle G(z), z\rangle \leq 0.
\]
Such an inequality translates in terms of expansion as
\begin{equation}\label{basic}
-|z_1|^2-|z_2|^2+\sum_{|\al|\geq 2} \Re q^1_\al z_1^{\al_1}\overline{z_1}z_2^{\al_2}+\sum_{|\al|\geq 2} \Re q^2_\al z_1^{\al_1}z_2^{\al_2}\overline{z_2}\leq 0.
\end{equation}
Replacing $(z_1,z_2)$ by $(e^{i\theta k_1}z_1, e^{i\theta k_2}z_2)$, with $\theta\in \R$ and $k_1,k_2\in \mathbb Z$, we obtain the expression
\begin{equation}\label{basic2}
\begin{split} 
-|z_1|^2-|z_2|^2+&\sum_{|\al|\geq 2, (\al_1-1)k_1+\al_2k_2=0} \Re q^1_\al z_1^{\al_1}\overline{z_1}z_2^{\al_2}\\&+\sum_{|\al|\geq 2, \al_1k_1+(\al_2-1)k_2=0} \Re q^2_\al z_1^{\al_1}z_2^{\al_2}\overline{z_2}+R(e^{i\theta})\leq 0,
\end{split}
\end{equation}
where $R(e^{i\theta})$ are harmonic terms with some common period. Integrating
\eqref{basic2} in $\theta$ over such a period causes the term $R(e^{i\theta})$ to
disappear, and we get a new expression 
\begin{equation}\label{basic3} 
\begin{split}
-|z_1|^2-|z_2|^2+&\sum_{|\al|\geq 2, (\al_1-1)k_1+\al_2k_2=0} \Re q^1_\al
z_1^{\al_1} {\overline{z_1}}z_2^{\al_2}\\&+\sum_{|\al|\geq 2, \al_1k_1+(\al_2-1)k_2=0} \Re q^2_\al z_1^{\al_1}z_2^{\al_2}{\overline{z_2}}\leq 0.
\end{split}
\end{equation}
This means that the vector field 
\[
G^{(k_1,k_2)}(z)=(-z_1+\sum_{|\al|\geq 2, (\al_1-1)k_1+\al_2k_2=0} q^1_\al z^\al, -z_2+\sum_{|\al|\geq 2, \al_1k_1+(\al_2-1)k_2=0} q^2_\al z^\al),
\]
is again a Herglotz vector field associated with the class $\mathcal M$.

\subsection{Slice reduction} Let $\|v\|=1$. Let $G(z)$ be an autonomous Herglotz vector field associated with the class $\mathcal M$.  For $\zeta\in \D$, let 
\[
-\zeta p_v(\zeta)=\langle G(\zeta v), v\rangle.
\]
It is easy to see that $p_v(\zeta)=1+\tilde{p}_v(\zeta)$ belongs to the Carath\'eodory class;  in particular, (see,  e.g., \cite{P}), its coefficients are bounded by $2$. A direct computation gives
\[
p_v(\zeta)=1-\sum_{m=1}^\infty \left(\sum_{|\al|=m+1}(q_\al^1 \overline{v_1}+q_\al^2 \overline{v_2})v^{\al}\right)\zeta^{m}.
\]
In particular, for all $m\in \N$, $m\geq 1$,
\begin{equation}\label{stima}
\sup_{\|v\|=1}\left|\sum_{|\al|=m+1}(q_\al^1 \overline{v_1}+q_\al^2 \overline{v_2})v^{\al} \right|\leq 2.
\end{equation}
This condition is  necessary but not sufficient for $p_v$ to belong to the Carath\'eodory class. Observe that by \cite[Corollary 2.3]{P}, if for some $\|v\|=1$ and $m\geq 1$,
\[
\left|\sum_{|\al|=m+1}(q_\al^1 \overline{v_1}+q_\al^2 \overline{v_2})v^{\al} \right|= 2,
\]
then 
\[
p_v(\zeta)=\sum_{l=1}^m t_m \frac{e^{i\theta+2\pi i l/m}+z}{e^{i\theta+2\pi i l/m}-z}
\]
for some $\theta\in \R$ and $t_j\geq 0$ with $\sum t_j=1$.

A necessary and sufficient condition for $p_v$ to belong to the Carath\'eodory class is the following (\cite[Thm. 2.4]{P}). For all $m\geq 1$,
\begin{equation}\label{pom}
\sum_{k=0}^m \sum_{l=0}^m \left(\sum_{|\al|=k-l+1}(q_\al^1 \overline{v_1}+q_\al^2 \overline{v_2})v^{\al}\right)\lambda_k\overline{\lambda_l}\geq 0,
\end{equation}
for all $\lambda_0,\ldots, \lambda_m\in \C$, with the convention that for $k-l+1\leq -2$,
\[
\left(\sum_{|\al|=k-l+1}(q_\al^1 \overline{v_1}+q_\al^2 \overline{v_2})v^{\al}\right)=\overline{\left(\sum_{|\al|=l-k-1}(q_\al^1 \overline{v_1}+q_\al^2 \overline{v_2})v^{\al}\right)}
\] 
and $\sum_{|\al|\leq 1}(q_\al^1 \overline{v_1}+q_\al^2 \overline{v_2})v^{\al}=2$.

\section{Coefficient bounds:  $q^1_{m,0}$ and $b^1_{m,0}$} Let $G(z)$ be an autonomous Herglotz vector field associated with the class $\mathcal M$ and fix $m\in \N$, $m\geq 2$.  Using the trick of harmonic decoupling, consider $G^{(0,1)}(z)$. Then
\[
G^{(0,1)}(z_1,z_2)=(-z_1+\sum_{m\geq 2} q^1_{m,0}z_1^m, -z_2+\sum_{m\geq 2} q^2_{m,1}z_1^m z_2).
\]
From \eqref{stima} we obtain for all $m\geq 2$,
\begin{equation}\label{est1}
\sup_{\|v\|=1}\left|q_{m,0}^1 |v_1|^2+q_{m-1,1}^2|v_2|^2\right||v_1^{m-1}| \leq 2
\end{equation}
Taking $v_1=1, v_2=0$, we obtain 
\begin{equation}\label{est-m0}
|q^1_{m,0}|\leq 2.
\end{equation} 
This bound is sharp, as  can be seen by considering the autonomous Herglotz vector field $G(z_1,z_2)=(-z_1(1+z_1)(1-z_1)^{-1}, -z_2)$.

Now, for $m=2$, from \eqref{ODE}, we obtain
\[
a^1_{2,0}(t)=e^{-t}\int_0^t e^{-\tau}q^1_{2,0}(\tau)d\tau.
\]
By \eqref{est-m0} and \eqref{b-a}, we then have  for all $f\in S^0$,
\[
|b^1_{2,0}|\leq 2.
\]
The bound is sharp, as one  sees by considering the map $(k(z_1), z_2)$, where $k(z_1)$ is the Koebe function in $\D$.

A similar bound for $|b^1_{m,0}|$ is not known. 

\begin{question}
Is it true that $|b_{m,0}^1|\leq m$ for all $f\in S^0$ and $m\in \N$, $m\geq 3$? 
\end{question}
Note that if the bound is correct, it is sharp  as one sees from  the function $(k(z_1), z_2)$, where $k(z_1)$ is the Koebe function in $\D$.

\section{Coefficient bounds:  $q^1_{0,m}$ and $b^1_{0,m}$}

Let $G(z)$ be an autonomous Herglotz vector field associated with the class $\mathcal M$. Fix $m\in \N$, $m\geq 2$. Using the decoupling harmonic terms trick, consider the  vector field  $G^{(m,1)}$, given by
\[
G^{(m,1)}(z_1,z_2)=(-z_1+q^1_{0,m}(t)z_2^m, -z_2).
\]
Since $-G^{(m,1)}\in \mathcal M$, imposing the condition $\Re \langle G^{(m,1)}(z), z\rangle \leq 0$, we get 
\[
-|z_1|^2-|z_2|^2+\Re q^1_{0,m} \overline{z_1}z_2^m \leq 0.
\]
Setting $z_1=xe^{i(\theta+\eta)}$, $z_2=ye^{i\theta/{m}}$ with $x,y\geq 0$ and
$q^1_{0,m}e^{-i\eta}=|q^1_{0,{m}}|$,  we obtain the equivalent 
 equation
\[
-x^2-y^2+|q^1_{0,m}|xy^m \leq 0, \quad x,y\geq 0, x^2+y^2\leq 1.
\]
Using the method of Lagrange multipliers, one checks easily that the maximum
for the function $(x,y)\mapsto -x^2-y^2+|q^1_{0,m}|xy^m$ under the
constraint  $x,y\geq 0, x^2+y^2\leq 1$ is attained at the point $x=\frac{1}{\sqrt{(1+m)}}, y=\sqrt{\frac{m}{(1+m)}}$. Hence the previous inequality is satisfied if and only if
\begin{equation}
\label{est-0m}
|q^1_{0,m}|\leq \frac{(1+m)^{\frac{m+1}{2}}}{m^{\frac{m}{2}}}.
\end{equation}
Note that \eqref{est-0m} gives the sharp bound for the coefficients $q^1_{0,m}$ in the class $\mathcal M$.

As before, for $m=2$, from \eqref{ODE},  \eqref{est-m0} and \eqref{b-a} we then have for all $f\in S^0$,
\[
|b^1_{0,2}|\leq \frac{3\sqrt{3}}{2}.
\]
In particular, the map
\[
\Phi:(z_1,z_2)\mapsto (z_1+\frac{3\sqrt{3}}{2} z_2^2, z_2)\in S^0,
\]
is a {\sl bounded} support point. Note also that $\Phi$ is an automorphism of
$\C^2$ and provides an affirmative answer to Question \ref{questionauto} for this  automorphism.  

This result, together with the germinal idea of decoupling harmonic terms, was proved in \cite{B}.

\begin{question}
Is it true that 
\[
|b_{0,m}^1|\leq u_m:=\frac{(1+m)^{\frac{m+1}{2}}}{m^{\frac{m}{2}}}
\]
 for all $f\in S^0$ and $m\in \N$, $m\geq 3$?
\end{question}
 Note that if the bound is correct, it is sharp  (consider the function $(z_1+u_mz_2^m, z_2)$).

\bibliographystyle{amsplain}

\end{document}